# Derivatives and integrals: matrix order operators as an extension of the fractional calculus

Cleiton Bittencourt da Porciúncula

UERGS, State University of Rio Grande do Sul, Brazil

Independência, Av., Santa Cruz do Sul Campus, ZIP Code: 96816-501

Corresponding author:

Phone/Fax Number: +55 51 3715-6926

Email: cleiton-porciuncula@uergs.edu.br

Abstract.

A natural consequence of the fractional calculus is its extension to a matrix order of differentiation and integration. A matrix-order derivative definition and a matrix-order integration arise from the generalization of the gamma function applied to the fractional differintegration definition. This work focuses on some results applied to the Riemann-Liouville version of the fractional calculus extended to its matrix-order concept. This extension also may apply to other versions of fractional calculus. Some examples of possible ordinary and partial matrix-order differential equations and related exact solutions are shown.

1. Introduction

There are several good works and reviews about fractional derivatives and integrals. In a general fashion, we may consider a general fractional derivative operator as:

$$_aD^\alpha f(x) = \frac{d^\alpha}{dx^\alpha}f(x) = \frac{d^p}{dx^p}\int_a^x K_D(x,t,\alpha)\phi_D(f(t),f^n(t))dt \qquad (1.1)$$

and a fractional integral operator as:

$$J^\alpha f(x) = {}_aD^{-\alpha}f(x) = \int_a^x K_J(x,t,\alpha)\phi_J(f(t),f^n(t))dt \qquad (1.2)$$

with $\alpha \in \mathbb{C}$, $t, x, a \in \mathfrak{R}$ and $p \in \mathbb{N}$. In eqs. (1.1) and (1.2), the kernels $K_D$ and $K_f$ vary according to the version of the fraction differintegration. Functions $\phi_D$ and $\phi_J$ also may assume any form as a function of $f(t)$ and the $n$-derivatives of $f(t)$, $f^n(t)$. Table 1 shows different forms of these operators under different kernels and shapes of integrable functions, with some corresponding references.

| $\alpha$-derivative | $\alpha$-integral | |
|---|---|---|
| $_cD_x{}^\alpha f(x) = \frac{1}{\Gamma(n-\alpha)}\frac{d^n}{dx^n}\int_c^x (x-t)^{(n-\alpha-1)} f(t)dt$; $x > c, n > \alpha + 1$ | $_cD_x{}^{-\alpha} f(x) = \frac{1}{\Gamma(\alpha)}\int_c^x (x-t)^{(\alpha-1)} f(t)dt$; $x > c$ | a |
| $_{-\infty}D_x{}^\alpha f(x) =$ $= \frac{1}{\Gamma(n-\alpha)}\frac{d^n}{dx^n}\int_{-\infty}^x (x-t)^{(n-\alpha-1)} f(t)dt$; $n > \alpha + 1$ | $_{-\infty}D_x{}^{-\alpha} f(x) = \frac{1}{\Gamma(\alpha)}\int_{-\infty}^x (x-t)^{(\alpha-1)} f(t)dt$ | b |
| $\frac{1}{\Gamma(n-\alpha)}\left(\frac{d}{dx}\right)^p \int_a^x (x-t)^{-(\alpha-n+1)} f(t)dt$; $x > a$; $n-1 < \alpha \le n$ $\frac{1}{\Gamma(n-\alpha)}\left(-\frac{d}{dx}\right)^p \int_x^b (x-t)^{-(\alpha-n+1)} f(t)dt$; $x < b$; $n-1 < \alpha \le n$ | $_0D_x{}^{-\alpha} f(x) = \frac{1}{\Gamma(\alpha)}\int_0^x (x-t)^{(\alpha-1)} f(t)dt$ | c |
| $D^\alpha f(x) = \frac{1}{2\pi}\int_0^{2\pi} f(x-t)\psi_{1-\alpha}'(t)dt$; $0 < \alpha < 1$ $\psi^\alpha(t) = \sum_{k=-\infty, k\ne 0}^{\infty} \frac{\exp(ikt)}{(\pm ik)^\alpha}$ | $_0D_x{}^{-\alpha} f(x) = \frac{1}{\Gamma(\alpha)}\int_x^\infty (x-t)^{(\alpha-1)} f(t)dt$ | d |
| $D_{a+}^\alpha f(x) = \frac{1}{\Gamma(n-\alpha)}\int_a^x (x-t)^{-(\alpha-n+1)} f^n(t)dt$; $x > a$; $n-1 < \alpha \le n$ $D_{b-}^\alpha f(x) = \frac{(-1)^n}{\Gamma(n-\alpha)}\int_x^b (x-t)^{-(\alpha-n+1)} f^n(t)dt$; $x < b$; $n-1 < \alpha \le n$ | $_aD_x{}^{-\alpha} f(x) = \frac{1}{\Gamma(\alpha)}\int_a^x (x-t)^{(\alpha-1)} f(t)dt$ | e |

Table 1 – Different versions of the fractional differintegral operators ($Re(\alpha > 0)$: a) Riemann: [1] Miller and Ross, 1993; b) Liouville: [1]Miller and Ross, 1993, [2] Oliveira e Machado, 2014; c) Riemann-Liouville, Abel equation: [1]Miller and Ross, 1993, [3]Trujillo et al, 2006, [3] [4]Samko, Kilbas and Marichev; d) Weyl: [5]Rudolf Hilfer, 1993, [6]Fausto Ferrari; e) Caputo: [2] Trujillo et al, 2006, [7] Caputo, 1966, [8] Caputo, 1967, [9]Kumar et al, 2019.

Starting from the Rieman-Liouville version by assuming $n = 1$, $p = 1$, with the lower integration limit $a = 0$, the fractional derivative becomes:

$$D^\alpha f(x) = \frac{d^\alpha}{dx^\alpha} f(x) = \frac{1}{\Gamma(1-\alpha)}\frac{d}{dx}\int_0^x (x-t)^{-\alpha} f(t)dt \qquad (1.3)$$

The corresponding fractional integral is:

$$D^{-\alpha}f(x) = J^{\alpha}f(x) = \frac{1}{\Gamma(\alpha)}\int_0^x (x-t)^{(\alpha-1)}f(t)dt \tag{1.4}$$

In this work, the following developments base on eqs. (1.3), (1.4). A natural extension of definitions (1.1)-(1.4) is possible if thinking in the sense of a differentiation/integration order as a squared matrix as follows:

$$\frac{d^M}{dx^M}f(x), \text{ with } M = \begin{bmatrix} a_{11} & a_{12} & \cdots & a_{1n} \\ a_{21} & a_{22} & \cdots & a_{2n} \\ \vdots & \vdots & \vdots & \vdots \\ a_{n1} & a_{n2} & \cdots & a_{nn} \end{bmatrix}$$

$M$ is a matrix $n \times n$ with, $n \in \mathbb{Z}$ and $a_{ij,\ 1<i<n,1<j<n}$ elements; $a \in \mathbb{C}$, and $f(x)$ a continuous scalar or matrix real or complex function function ($x \in \mathbb{C}$). The "matrix-order differentiation" may look improbable at a first glance; however, it is possible to define it also in terms of the gamma function by simply following the known definitions of fractional derivative and integral in eqs. (1.3), (1.4). The differintegration operators expressed as a matrix order (in terms of a generic kernel) become:

$$_aD^M F(x) = \frac{d^M}{dx^M}F(x) = \frac{d^p}{dx^p}\int_a^x K_D(x,t,M)\phi_D(F(t),F^n(t))dt \tag{1.5}$$

$$J^M f(x) = {}_aD^{-M}f(x) = \int_a^x K_J(x,t,M)\phi_J(F(t),F^n(t))dt \tag{1.6}$$

Expressing the matrix-order differintegration operators as a Rieman-Liouville according to eqs. (1.3) - (1.4), one obtains:

$$D^M f(x) = \frac{d^M}{dx^M}f(x) = [\Gamma^{-1}(I-M)]\frac{d}{dx}\int_0^x (x-t)^{-M}F(t)dt \tag{1.7}$$

$$D^{-M}f(x) = J^M f(x) = [\Gamma^{-1}(M)]\int_0^x (x-t)^{(M-I)}F(t)dt \tag{1.8}$$

where $I$ is the identity matrix and $M \in \mathbb{C}^{n \times n}$. In eqs. (1.5) - (1-8), the function $f$ is now expressed with an uppercase $F$ to account the possibility that the function $F$ may be a matrix of functions (since that $\dim F = \dim M$), that is, the number of rows and columns in $F$ is equal to those in $M$, both with the same dimension), and not only a single function $f$. $F$ may also be a column vector of functions, provided that the number of rows of $F$ has the same number of columns of $M$ or $(x - t)^{-M}$. This extension of the fractional calculus to a matrix order calculus is not a new idea. Naber develops this theme in their works [10], [11]. Phillips [12], [13] uses the concept of matrix order differentiation and integration in the development of the econometric theory and advanced probability distributions. However, there are no references, for example, about developments of exact or numeric matrix-order differential equations, other related problems of such objects, and possible applications. Before continuing the analysis of this matrix order calculus, Table 2 summarizes the nomenclature and symbols used along with this text.

| Notation / Symbol | Used for: |
|---|---|
| $M$ | The uppercase letter indicates a square matrix with $n \times n$ elements. |
| $I$ | Identity matrix with $n \times n$ elements. |
| $f(x), f(t)$ | A single function depending on $x$ or $t$ |
| $F(x), F(t)$ | A matrix of functions $f_{ij}(x)$ or $f_{ij}(t)$ in the form: $F = \begin{bmatrix} f_{11} & \cdots & f_{1n} \\ \vdots & \ddots & \vdots \\ f_{n1} & \cdots & f_{nn} \end{bmatrix}$ or a column vector function in the form: $F = \begin{bmatrix} f_1 \\ \vdots \\ f_m \end{bmatrix}$ where $n = m$, i.e., the number of columns of $M$ is equal to the number of rows of $F$. |
| $D^M, \frac{d^M}{dx^M}$ | The $M$-derivative of a specified function. |
| $\partial^M, \frac{\partial^M}{\partial x^M}, \partial_x^M$ | The partial $M$-derivative of a specified function concerning some variable $x$. |
| $D^{-M}, J^M, \int dx^M$ | The $M$-integral of a specified function |
| $\Gamma(A)$ | The matrix gamma function of some square matrix $A$. |
| $B(A)$ | The matrix beta function of some square matrix $A$. |
| $\Gamma^{-1}(A)$ | The inverse of the matrix gamma function of some square matrix $A$. |
| $B^{-1}(A)$ | The inverse of the matrix beta function of some square matrix $A$. |
| $K^{-A}(x,t)$ | Equal to $[(K(x,t))^A]^{-1}$ or $inv\,[(K(x,t)^A]$ |

Table 2 – Notation used in this article.

In eqs. (1.7) and (1.8), the matrix gamma functions are:

$$\Gamma(I - M) = \int_0^{+\infty} t^{-M} \exp(-t)\, dt \qquad (1.9)$$

$$\Gamma(M) = \int_0^{+\infty} t^{(M-I)} \exp(-t)\, dt \qquad (1.10)$$

Eqs. (1.9) – (1.10) are the matrix gamma function, and its properties are broadly covered in the literature, as well as several properties and applications derived from those definitions, especially applied to multivariate statistics [14], [15], [16], [17].

The term $t^{-M}$ in eq. (1.9) is:

$$t^{-M} = \exp\left(\ln t^{-M}\right) = \exp(-M \ln(t)) \qquad (1.11)$$

The exponential of a matrix is a well-known result based on the decomposition of eigenvalues and eigenvectors of a diagonalizable matrix [19]:

$$\exp(At') = V \exp(\Lambda t')\, V^{-1} \qquad (1.12)$$

with $t' = \ln t$. For a $n \times n$ matrix, one obtains:

$$\exp(\Lambda t') = \begin{bmatrix} \exp(\lambda_1 t') & \cdots & 0 \\ \vdots & \ddots & \vdots \\ 0 & \cdots & \exp(\lambda_n t') \end{bmatrix} \qquad (1.13)$$

The corresponding eigenvector matrix is:

$$V = \begin{bmatrix} v_1 & \cdots & v_n \end{bmatrix} = \begin{bmatrix} v_{11} & \cdots & v_{1n} \\ \vdots & \ddots & \vdots \\ v_{n1} & \cdots & v_{nn} \end{bmatrix} \qquad (1.14)$$

The vectors $v_1$ up to $v_n$ are the associated eigenvectors with the eigenvalues $\lambda_1$ up to $\lambda_n$. These eigenvalues are computed as:

$$\det(\lambda I - A) = 0 \qquad (1.15)$$

The eigenvalues $\lambda_i$ associated with each eigenvector $v_i$ is determined via the well-known relation $Av_i = \lambda_i v_i$, being for eq. (1.11), $A = -M$, and changing (1.15) to:

$$\det(\lambda I + M) = 0 \qquad (1.16)$$

When applying the decomposition of matrices to eq. (1.11), the result is:

$$\exp(-M \ln t) = V \exp(\Lambda t') V^{-1} \qquad (1.17)$$

Now with $V$ as the matrix of eigenvectors of $\exp(At')$, according to (1.14):

$$\exp(\Lambda t') = \exp(-M \ln t) =$$
$$= \begin{bmatrix} \exp(\lambda_1 \ln t) & \cdots & 0 \\ \vdots & \ddots & \vdots \\ 0 & \cdots & \exp(\lambda_n \ln t) \end{bmatrix} = \begin{bmatrix} t^{\lambda_1} & \cdots & 0 \\ \vdots & \ddots & \vdots \\ 0 & \cdots & t^{\lambda_n} \end{bmatrix} \qquad (1.18)$$

The same procedure applies to $t^{M-I}$ in eq. (1.10):

$$t^{(M-I)} = \exp(\ln t^{(M-I)}) = \exp\big((M - I) \ln(t)\big) \qquad (1.19)$$

$$\exp(\Upsilon t') = \exp\big((M - I) \ln(t)\big) =$$
$$= \begin{bmatrix} \exp(\upsilon_1 \ln t) & \cdots & 0 \\ \vdots & \ddots & \vdots \\ 0 & \cdots & \exp(\upsilon_n \ln t) \end{bmatrix} = \begin{bmatrix} t^{\upsilon_1} & \cdots & 0 \\ \vdots & \ddots & \vdots \\ 0 & \cdots & t^{\upsilon_n} \end{bmatrix} \qquad (1.20)$$

but now, $A = M - I$ and:

$$\exp(\Upsilon t') = \exp\big((M - I) \ln(t)\big) = W \exp(\Upsilon t') W^{-1} \qquad (1.21)$$

Now with $W$ as the matrix of eigenvectors of $\exp(\Upsilon t')$:

$$W = [w_1 \ \cdots \ w_n] = \begin{bmatrix} w_{11} & \cdots & w_{1n} \\ \vdots & \ddots & \vdots \\ w_{n1} & \cdots & w_{nn} \end{bmatrix} \tag{1.22}$$

The eigenvalues $v_i$ associated with each eigenvector $w_i$ is determined via the known relation $Aw_i = v_i w_i$, being in this case, $A = M - I$. The eigenvalues are computed as eq. (1.12):

$$\det(vI - (M - I)) = \det(I(v + 1) - M) = 0 \tag{1.23}$$

For the term $(x - t)^{-M}$ in eq. (1.7), working first with $(x - t)^M$ the same procedure takes to:

$$(x - t)^M = \exp(\ln(x - t)^M) = \exp(M \ln(x - t)) \tag{1.24}$$

and

$$\det(hI - M) = 0 \tag{1.25}$$

with the following diagonalizable matrix:

$$\exp(Ht^*) = \exp(M \ln(x - t)) = \begin{bmatrix} e^{h_1 \ln(x-t)} & \cdots & 0 \\ \vdots & \ddots & \vdots \\ 0 & \cdots & e^{h_n \ln(x-t)} \end{bmatrix} =$$

$$= \begin{bmatrix} (x - t)^{h_1} & \cdots & 0 \\ \vdots & \ddots & \vdots \\ 0 & \cdots & (x-)t^{h_n} \end{bmatrix} \tag{1.26}$$

Now with $P$ as the matrix of eigenvectors of $\exp(Ht^*)$ and doing as $t^* = \ln(x - t)$:

$$P = [p_1 \ \cdots \ p_n] = \begin{bmatrix} p_{11} & \cdots & p_{1n} \\ \vdots & \ddots & \vdots \\ p_{n1} & \cdots & p_{nn} \end{bmatrix} \tag{1.27}$$

$$\exp(M \ln(x - t)) = P \exp(Ht^*) P^{-1} \qquad (1.28)$$

The eigenvalues $h_i$ associated with each eigenvector $p_i$ is again determined by using the relation $Ap_i = h_i p_i$, which leads to eq. (1.25), being in this case, $A = M$.

For the term $(x - t)^{(M-I)}$ in eq. (1.8), the same previous procedures take to:

$$(x - t)^{(M-I)} = \exp\bigl(\ln(x - t)^{(M-I)}\bigr) = \exp\bigl((M - I) \ln(x - t)\bigr) \qquad (1.29)$$

The corresponding eigenvalues are computed as previously done for eq. (1.23):

$$\det\bigl(\omega I - (M - I)\bigr) = \det(I(\omega + 1) - M) = 0 \qquad (1.30)$$

The corresponding diagonalizable matrices are:

$$\exp(\Omega t') = \exp\bigl((M - I) \ln(x - t)\bigr) = \begin{bmatrix} e^{\omega_1 \ln(x-t)} & \cdots & 0 \\ \vdots & \ddots & \vdots \\ 0 & \cdots & e^{\omega_n \ln(x-t)} \end{bmatrix}$$

$$= \begin{bmatrix} (x - t)^{\omega_1} & \cdots & 0 \\ \vdots & \ddots & \vdots \\ 0 & \cdots & (x-)t^{\omega_n} \end{bmatrix} \qquad (1.31)$$

Now with $Q$ as the matrix of eigenvectors of $\exp(\Omega t')$:

$$Q = [q_1 \ \cdots \ q_n] = \begin{bmatrix} q_{11} & \cdots & q_{1n} \\ \vdots & \ddots & \vdots \\ q_{n1} & \cdots & q_{nn} \end{bmatrix} \qquad (1.32)$$

$$\exp\bigl((M - I) \ln(x - t)\bigr) = Q \exp(\Omega t') Q^{-1} \qquad (1.33)$$

The eigenvalues $\omega_i$ associated with each eigenvector $q_i$ is also determined by using the relation $Aq_i = \omega_i q_i$, being in this case, $A = M - I$.

So, combining eqs. (1.7), (1.9):

$$D^M f(x) = \frac{d^M}{dx^M} f(x) =$$

$$= \left[ \int_0^{+\infty} \begin{bmatrix} v_{11} & \cdots & v_{1n} \\ \vdots & \ddots & \vdots \\ v_{n1} & \cdots & v_{nn} \end{bmatrix} \begin{bmatrix} t^{\lambda_1} & \cdots & 0 \\ \vdots & \ddots & \vdots \\ 0 & \cdots & t^{\lambda_n} \end{bmatrix} \begin{bmatrix} v_{11} & \cdots & v_{1n} \\ \vdots & \ddots & \vdots \\ v_{n1} & \cdots & v_{nn} \end{bmatrix}^{-1} e^{-t} dt \right]^{-1}$$

$$\frac{d}{dx} \int_0^x \left\{ \begin{bmatrix} p_{11} & \cdots & p_{1n} \\ \vdots & \ddots & \vdots \\ p_{n1} & \cdots & p_{nn} \end{bmatrix} \begin{bmatrix} (x-t)^{h_1} & \cdots & 0 \\ \vdots & \ddots & \vdots \\ 0 & \cdots & (x-)t^{h_n} \end{bmatrix} \begin{bmatrix} p_{11} & \cdots & p_{1n} \\ \vdots & \ddots & \vdots \\ p_{n1} & \cdots & p_{nn} \end{bmatrix}^{-1} \right\}^{-1} F(t) dt$$

(1.34)

Now, combining (1.8) and (1.10) for the matrix-order integral:

$$J^M f(x) = \int f(x) dx^M =$$

$$= \left[ \int_0^{+\infty} \begin{bmatrix} w_{11} & \cdots & w_{1n} \\ \vdots & \ddots & \vdots \\ w_{n1} & \cdots & w_{nn} \end{bmatrix} \begin{bmatrix} t^{v_1} & \cdots & 0 \\ \vdots & \ddots & \vdots \\ 0 & \cdots & t^{v_n} \end{bmatrix} \begin{bmatrix} w_{11} & \cdots & w_{1n} \\ \vdots & \ddots & \vdots \\ w_{n1} & \cdots & w_{nn} \end{bmatrix}^{-1} e^{-t} dt \right]^{-1}$$

$$\left[ \int_0^x \begin{bmatrix} q_{11} & \cdots & q_{1n} \\ \vdots & \ddots & \vdots \\ q_{n1} & \cdots & q_{nn} \end{bmatrix} \begin{bmatrix} (x-t)^{\omega_1} & \cdots & 0 \\ \vdots & \ddots & \vdots \\ 0 & \cdots & (x-)t^{\omega_n} \end{bmatrix} \begin{bmatrix} q_{11} & \cdots & q_{1n} \\ \vdots & \ddots & \vdots \\ q_{n1} & \cdots & q_{nn} \end{bmatrix}^{-1} F(t) dt \right] \quad (1.35)$$

Eqs. (1.34) and (1.35) are complex expressions difficult to determine analytically even for simple functions and for small matrices. Numerical approaches are a good way to find approximations for these expressions.

2. Some basic properties of the matrix differintegration operators

2.1 Additive property

From the fractional calculus [1], [5], we know that:

$$D^{-\alpha}[D^{-\beta} f(t)] = D^{-(\alpha+\beta)} f(t) = D^{-\beta}[D^{-\alpha} f(t)] \quad (\alpha, \beta > 0) \tag{2.1}$$

Or, using alternative notation (Table 2):

$$J^\alpha[J^\beta f(t)] = J^{(\alpha+\beta)} f(t) = J^\beta[J^\alpha f(t)] \quad (\alpha, \beta > 0) \tag{2.2}.$$

Except for the fact that, using notation $D$ as expressed in eq. (2.1), the minus sign indicates the opposite operation to differentiation, that is, integration, both notations mean that: 1) the order of fractional integration leads to the same result, and 2), multiple fractional integrations are the final integration order corresponding to the summation of each fractional integral order. This is also valid if we consider matrix-integration, by generalizing the proof for (2.1) or (2.2) as follows:

*Proof*:

$$D^{-M}D^{-N}F(x) = \Gamma^{-1}(M) \int_a^x (x-u)^{(M-I)} [D^{-N}F(u)] du \qquad (2.3)$$

$$D^{-M}D^{-N}F(x) = \Gamma^{-1}(M) \int_a^x (x-u)^{(M-I)} \left[\Gamma^{-1}(N) \int_a^u (u-t)^{(N-I)} dt\, F(u)\right] du\, dt \qquad (2.4)$$

$$D^{-M}D^{-N}F(x) = \Gamma^{-1}(M)\Gamma^{-1}(N) \int_a^x (x-u)^{(M-I)} \left[\int_t^x (u-t)^{(N-I)} F(t) dt\right] du \qquad (2.5)$$

$$D^{-M}D^{-N}F(x) = \Gamma^{-1}(M)\Gamma^{-1}(N) \int_a^x F(t) \left[\int_t^x (x-u)^{(M-I)} (u-t)^{(N-I)} dt\right] du \qquad (2.6)$$

Making the substitution $= \frac{u-t}{x-t}$:

$$D^{-M}D^{-N}F(x) = \Gamma^{-1}(M+N) \int_a^x (x-t)^{(M+N-I)} F(t) \left[\int_0^1 (1-y)^{(M-I)} y^{(N-I)} dy\right] du$$

(2.7)

The inner integral is the beta function in the form of the matrix beta function:

$$\mathrm{B}(M,N) = \int_0^1 (1-y)^{(M-I)} y^{(N-I)} dy \qquad (2.8)$$

The relation between gamma and beta function in matrix form is similar to scalar form [14]:

$$\mathrm{B}(M,N) = \Gamma(M)\Gamma(N)\Gamma^{-1}(M+N) \qquad (2.9)$$

It is important to emphasize, that the matrices $M$ and $N$ in $C^{n \times n}$ must be diagonalizable. Also, $MN = NM$, that is, the matrices must be commutable and all their eigenvalues must present real part greater than zero, $Re(\lambda_M, \lambda_N) > 0$, being $\lambda$ the corresponding eigenvalues for $M$ and $N$.

$$D^{-M}D^{-N}F(x) = \Gamma^{-1}(M+N)\left[\int_a^x (x-t)^{(M+N-I)}F(t)\right]\text{B}(M,N) \qquad (2.10)$$

Replacing (2.9) into (2.10):

$$D^{-M}D^{-N}F(x) = \Gamma^{-1}(M+N)\Gamma(M)\Gamma(N)\Gamma^{-1}(M+N)\left[\int_a^x (x-t)^{(M+N-I)}F(t)\right] =$$

$$= \Gamma(M)\Gamma(N)\left[\int_a^x (x-t)^{(M+N-I)}F(t)\right] = D^{-(M+N)}$$

$$(2.11)$$

Therefore:

$$D^{-M}D^{-N}F(x) = D^{-(M+N)}F(x) \qquad (2.12)$$

Which completes the proof.

Unfortunately, the additive property for matrix-integral operators is not completely valid for matrix differential operators, because it is not valid for their analogous fractional operators [1], as observed in [3]:

$$D_{a+}^{\alpha}D_{a+}^{\beta}f(x) = D_{a+}^{\alpha+\beta}f(x) - \sum_{j=1}^{m}(D_{a+}^{\beta-j}f)(a+)\frac{(x-\alpha)^{-j-\alpha}}{\Gamma(1-j-\alpha)} \qquad (2.13)$$

With, $\beta > 0$, $-1 < \alpha \leq n$, $m - 1 < \beta \leq m$, $n, m \in \mathbb{N}$, and $\alpha + \beta \leq n$.

2.2 Inverse operation

In the classical integer order calculus, we know that $\frac{d}{dx}\int f(x)dx = f(x)$, that is, the derivative of the integral of a function returns the original function. This is also valid for the fractional calculus. The fraction differentiation is the opposite of the fraction integration [4] so that the $\alpha$-derivative of the $\alpha$-integral of $f(x)$ returns $f(x)$ [18]:

$$D^\alpha D^{-\alpha} f(x) = f(x) \tag{2.14}$$

In terms of matrix differintegration (2.14) is:

$$D^M D^{-M} F(x) = D^M J^M F(x) = F(x) \tag{2.15}$$

*Proof:*

Writing (2.15) based on the definitions (1.7) and (1.8):

$$D^M D^{-M} F(x) = \Gamma^{-1}(I-M) \frac{d}{dx} \int_0^x (x-u)^{-M} \, \Gamma^{-1}(M) \int_0^u (u-t)^{(M-I)} F(u) \, du \, dt \tag{2.16}$$

Interchanging the integration order:

$$D^M D^{-M} F(x) = \Gamma^{-1}(I-M) \Gamma^{-1}(M) \frac{d}{dx} \int_0^x F(t) \int_t^x (x-u)^{-M} (u-t)^{(M-I)} \, dt \, du \tag{2.17}$$

Making the same substitution $= \frac{u-t}{x-t}$:

$$D^M D^{-M} F(x) = \Gamma^{-1}(I-M) \Gamma^{-1}(M) \frac{d}{dx} \int_0^x F(t) \int_0^1 (1-y)^{-M} y^{(M-I)} \, dy \, du \tag{2.18}$$

Applying the definition of the matrix beta function:

$$\mathrm{B}(M,N) = \mathrm{B}(M, I-M) = \int_0^1 y^{(M-I)} (1-y)^{(I-M-I)} dt = \int_0^1 y^{(M-I)} (1-y)^{-M} dt \tag{2.19}$$

Furthermore, based on eq. (2.9):

$$\mathrm{B}(M,N) = \mathrm{B}(M, I-M) = \Gamma(M)\Gamma(I-M)\Gamma^{-1}(M+I-M) = \Gamma(M)\Gamma(I-M)\Gamma^{-1}(I) =$$
$$= \Gamma(M)\Gamma(I-M) \tag{2.20}$$

Replacing (2.20) and (2.19) into eq. (2.18) one obtains:

$$D^M D^{-M} F(x) = \Gamma^{-1}(I-M)\Gamma^{-1}(M)\frac{d}{dx}\int_0^x F(t)\,\mathrm{B}\,du =$$

$$= \Gamma^{-1}(I-M)\Gamma^{-1}(M)\frac{d}{dx}\int_0^x F(t)\,\Gamma(M)\Gamma(I-M) \tag{2.21}$$

By canceling out the beta function with the inverse gamma functions, we have:

$$D^M D^{-M} F(x) = \frac{d}{dx}\int_0^x F(t)\,dt = F(x) \tag{2.22}$$

Which completes the proof.

2.3 Inverse operation with different matrix order of derivatives and integrals

Now let us consider the following property described for both cases, $q \geq p \geq 0$ and $p > q \geq 0$) [3], [18]:

$$D^p D^{-q} f(t) = D^{p-q} f(t) \tag{2.23}$$

The matrix form of the property (2.23) follows:

$$D^M D^{-N} F(x) = D^M \left( D^{-M} D^{-(N-M)} F(x) \right), M \neq N \tag{2.24}$$

First proof:
Following property (2.15) and applying to (2.24):

$$D^M D^{-N} F(x) = D^{(M-M)} \left( D^{(M-N)} F(x) \right) = D^{(M-N)} F(x) \tag{2.25}$$

An alternative proof is similar to the previous proof on $D^M D^{-M}$ and $D^{-M} D^{-N}$, if considering that the differentiation order $M$ is different from the integration order $N$:

Second proof:
Applying definitions (1.7), (1.8):

$$D^M D^{-N} F(x) = \Gamma^{-1}(I-(M-N))\frac{d}{dx}\int_0^x (x-t)^{-(N-M)} = D^M D^{-N} F(x) =$$

$$= \Gamma^{-1}(I - M + N) \frac{d}{dx} \int_0^x (x - t)^{-(N-M)} \tag{2.26}$$

$$D^M D^{-N} F(x) = \Gamma^{-1}(I - M)\Gamma^{-1}(N) \frac{d}{dx} \int_0^x (x - u)^{-M} \int_0^u (u - t)^{(N-I)} F(u) du dt \tag{2.27}$$

Replacing the differentials and change the variable $y = \frac{u-t}{x-t}$, we have:

$$D^M D^{-N} F(x) = \Gamma^{-1}(I - M)\Gamma^{-1}(N) \frac{d}{dx} \int_0^x F(t)(x - u)^{-M} \int_t^x (u - t)^{(N-I)} dt du \tag{2.28}$$

$$D^M D^{-N} F(x) = \Gamma^{-1}(I - M)\Gamma^{-1}(N) \frac{d}{dx} \int_0^x F(t) \int_0^1 (x - t)^{(N-M)}(1 - y)^{-M} y^{(N-I)} dy du$$
$$\tag{2.29}$$

Again, by applying the new matrix Beta function definition:

$$B(N, I - M) = \int_0^1 y^{(N-I)}(1 - y)^{(I-M-I)} dy = \int_0^1 y^{(N-I)}(1 - y)^{-M} dy \tag{2.30}$$

Moreover:

$$B(N, I - M) = \Gamma(N)\Gamma(I - M)\Gamma^{-1}(N + I - M) \tag{2.31}$$

Replacing (2.31) and (2.30) in (2.29):

$$D^M D^{-N} F(x) = \Gamma^{-1}(I - (M - N)) \frac{d}{dx} \int_0^x F(t)(x - t)^{-(M-N)} dt =$$
$$= \Gamma^{-1}(I - M + N)) \frac{d}{dx} \int_0^x F(t)(x - t)^{(N-M)} dt = D^{(M-N)} \tag{2.32}$$

Which completes the proof.

3.1 Matrix order ordinary differential equations (MO-ODE)

3.1.1. Non-homogeneous equations

Consider the following equation ($\phi(x) \neq 0$):

$$\frac{d^M}{dx^M}F(x) = \phi(x) \qquad (3.1)$$

which also may be written as:

$$D^M F(x) = \phi(x) \qquad (3.2)$$

or simply as:

$$D^M F = \phi \qquad (3.3)$$

is a linear ordinary matrix order differential equation (LOMODE). If using the Rieman-Liouville definition, the known function $\phi(x)$ must be different of zero. Let us consider that $F$ depends of a function $G$ as follows:

$$F(x) = D^{-M}G(x) \qquad (3.4)$$

Substituting (3.4) into (3.1 – 3.3), and applying the property in eq. (2.22):

$$D^M(D^{-M}G(x)) = \phi(x) = G(x) \qquad (3.5)$$

Returning to (3.4) based on (3.5), and applying the M-integral definition (1.7):

$$F(x) = D^{-M}\phi(x) \qquad (3.6)$$

$$F(x) = \Gamma^{-1}(M) \int_0^x (x-t)^{(M-I)} \phi(t)\, dt \qquad (3.7)$$

Now let us consider the following equation:

$$C_m D^M F(x) + C_n D^N F(x) = \phi(x) \qquad (3.8)$$

with $C_m$ and $C_n$ are a matrix of constant parameters with the same dimension as $M$ and $N$ ; $\dim(C_m, C_n) = \dim(M, N)$. $C_m$ and $C_n$ may also be a single constant parameter, so that $C_m = c_m$ and $C_N = c_n$. Eq. (3.8) is a non-homogeneous linear matrix order differential equation (NH-LMODE). In a similar fashion to solve (3.1), let us consider the same assumption in eq. (3.4): $F(x) = D^{-M}G(x)$ into eq. (3.8):

$$C_m D^M D^{-M} G(x) + C_n D^N D^{-M} G(x) = \phi(x) \qquad (3.9)$$

Applying the property in eq. (2.15):
$$C_m G(x) + C_n D^N D^{-M} G(x) = \phi(x) \qquad (3.10)$$

Eq. (3.10) presents a solution by two ways. The first one is considering that matrix $M$ is the composition of other two matrices, making $M = A + N$. The second form uses property in eq. (2.25), and solves the resulting equation $C_m G(x) + C_n D^{(N-M)} G(x) = \phi(x)$. The disadvantage of this second method is the presence of a second-order matrix ordinary differential equation. The first method exposed as follows, generates a linear first order matrix differential equation, most simple to solve. Then, replacing $M$ by $N + A$:

$$C_m G(x) + C_n D^N D^{-A-N} G(x) = \phi(x) \qquad (3.11)$$

Following property (2.2), we know that $D^{-A} D^{-N} = D^{-N} D^{-A}$, and with property (2.15), $D^M D^{-M} F = F$:
$$C_m G(x) + C_n D^{-A} G(x) = \phi(x) \qquad (3.12)$$

From this point one, we also have two options: first, insert $\phi(x)$ and $C_m G(x)$ as derivatives inside the integral $D^{-A}$ and solve the resulting ordinary differential equation as $dG/dt$ that arises. The second option is to apply $D^A$ to both sides in eq. (3.12) and then insert the remaining terms inside the integral $D^A$. This second way is more complicated because the existence of the term $d/dx$ in the definition of $D^A$. Then, using the first method and opening eq. (3.12) with definition (1.8):

$$C_m G(x) + C_n \Gamma^{-1}(A) \int_0^x (x-t)^{(A-I)} G(t) = \phi(x) \qquad (3.13)$$

Inserting all the terms inside the integral:

$$\int_0^x \left[ C_n (x-t)^{(A-I)} \Gamma^{-1}(A) G(t) + C_m \frac{dG(t)}{dt} - \frac{d\phi(t)}{dt} \right] dt = 0 \qquad (3.14)$$

Provided that the functions $G(0) = 0$ and $\phi(0) = 0$, because $\int_0^x \frac{dG(t)}{dt} dt = G(t) - G(0)$, and the same valid when integrating $\phi$, to not arise additional terms in (3.14). If the integral in (3.14) is zero, then:

$$C_n (x-t)^{(A-I)} \Gamma^{-1}(A) G(t) + C_m \frac{dG(t)}{dt} - \frac{d\phi(t)}{dt} = 0 \qquad (3.15)$$

Eq. (3.15) is a linear ordinary differential equation with $\frac{d\phi(t)}{dt} \neq 0$.

$$\frac{dy}{dt} + p(t) y = q(t) \qquad (3.16)$$

With the solution:

$$y(t) = \exp\left(-\int p(t) dt\right) \left[ \int q(t) \left( \exp \int p(t) dt \right) dt + K \right] \qquad (3.17)$$

Rearranging (3.15) to become similar to (3.16):

$$\frac{dG(t)}{dt} + C_m^{-1} C_n (x-t)^{(A-I)} \Gamma^{-1}(A) G(t) = C_m^{-1} \frac{d\phi(t)}{dt} \qquad (3.18)$$

The solution of eq. (3.18) is:

$$G(t) = \exp\left(-\int C_m^{-1} C_n (x-t)^{(A-I)} \Gamma^{-1}(A) dt\right) \left( \int C_m^{-1} \frac{d\phi(t)}{dt} \exp\left(\int C_m^{-1} C_n (x-t)^{(A-I)} \Gamma^{-1}(A) dt\right) dt + K \right) \; (3.19)$$

with $K$ a matrix of constants and $A = M - N$. As initially, we want to know $F(x) = D^{-M}G(x)$:

$$F(x) = \Gamma^{-1}(M) \int_0^x (x-t)^{(M-I)} G(t) dt \qquad (3.20)$$

Eq. (3.20) is the solution of eq. (3.8)

Let us generalize eq. (3.8) in the following way:

$$\sum_{n=1}^{n} C_{m_n} D^{M_n} F(x) = \phi(x) \qquad (3.21a)$$

$$(C_{m_1} D^{M_1} + C_{m_2} D^{M_2} + \cdots + C_{m_n} D^{M_n}) F(x) = \phi(x) \qquad (3.21b)$$

Eq. (3.21) is a non-homogenous linear matrix order differential equation with different $M_i$ orders and with $\phi(x) \neq 0$. If applying the same previous steps for the case with two terms, we can express the matrix order $M_1$ in terms of the subsequent orders and ancillary matrices $A$:

$$\begin{aligned} M_1 &= M_2 + A_1 \\ M_1 &= M_3 + A_2 \\ M_1 &= M_4 + A_3 \\ &\vdots \\ M_1 &= M_n + A_{n-1} \end{aligned} \qquad (3.22)$$

Let us also assume that there exists a solution $F(x)$ such that $F(x) = D^{-M_1}G(x)$:

$$(C_{m_1} D^{M_1} D^{-M_1} + C_{m_2} D^{M_2} D^{-M_1} + \cdots + C_{m_n} D^{M_n} D^{-M_1}) G(x) = \phi(x) \qquad (3.23)$$

Replacing (3.22) into (3.23):

$$(C_{m_1} + C_{m_2} D^{M_2} D^{-M_2 - A_1} + \cdots + C_{m_n} D^{M_n} D^{-M_n - A_{n-1}}) G(x) = \phi(x) \qquad (3.24)$$

Using the property (2.15) as in the first term:

$$(C_{m_1} + C_{m_2} D^{-A_1} + \cdots + C_{m_n} D^{-A_{n-1}}) G(x) = \phi(x) \qquad (3.25)$$

Which may be expressed as:

$$C_{m_1}G(x) + \sum_{n=2}^{n} C_{M_n} D^{-A_{n-1}} G(x) = \phi(x) \qquad (3.26)$$

With the definition (1.8):

$$C_{m_1}G(x) + \sum_{n=2}^{n} \int_0^x C_{m_n} \Gamma^{-1}(A_{n-1})(x-t)^{(A_{n-1}-I)} G(t)\, dt = \phi(x) \qquad (3.27)$$

Again, putting the all the terms inside the integral:

$$\int_0^x \left[ \sum_{n=2}^{n} \left( C_{m_n} \Gamma^{-1}(A_{n-1})(x-t)^{(A_{n-1}-I)} G(t) dt \right) + C_{m_1} \frac{dG(t)}{dt} - \frac{d\phi(t)}{dt} \right] dt = 0$$

$$(3.28)$$

Which takes to:

$$\frac{dG(t)}{dt} + C_{m_1}^{-1} \sum_{n=2}^{n} \left( C_{m_n} \Gamma^{-1}(A_{n-1})(x-t)^{(A_{n-1}-I)} G(t) dt \right) = C_{m_1}^{-1} \frac{d\phi(t)}{dt}$$

$$(3.29)$$

Eq. (3.29) is essentially the same linear matrix differential equation as (3.15), but now $p(t)$ has more terms. The solution is:

$$G(t) = \exp\left(-\int C_m^{-1} \sum_{n=2}^{n} C_n (x-t)^{(A_{n-1}-I)} \Gamma^{-1}(A_{n-1}) dt\right) \left(\int C_m^{-1} \frac{d\phi(t)}{dt} \exp\left(\int C_m^{-1} \sum_{n=2}^{n} C_n (x-t)^{(A_{n_1}-I)} \Gamma^{-1}(A_{n_1}) dt\right) dt + K\right) \qquad (3.30)$$

With $A_{n_1} = M_1 - M_n$ and $\frac{d\phi(t)}{dt} \neq 0$.

Finally, returning to the desired solution $F(x) = D^{-M_1} G(x)$

$$F(x) = \Gamma^{-1}(M_1) \int_0^x (x-t)^{(M_1-I)} G(t) dt \qquad (3.31)$$

Consider the following equation:

$$C_m D^M F(x) = F(x) \qquad (3.32)$$

Proceeding with the same previous summarized below:

1 – Assume $F(x) = D^{-M}G(x)$

2 – Based on the properties (2.12), (2.15) or (2.25), write the equation in terms of $G(x)$.

3 – Put the free terms with $G(x)$ and other functions inside the integral, since such functions are zero for $t = 0$, to not create additional terms after integration (that is: $\int_0^x \frac{dG(t)}{dt} dt = G(x) - G(0)$, $G(0) = 0$).

4 – Solve the ordinary differential equation to obtain $G(t)$ in the domain $t$ by applying indefinite integration.

5 – Solve $F(x) = D^{-M}G(x)$ by applying definition (1.8).

When applying the steps (1-4) in eq. (3.32), we have:

$$\Gamma^{-1}(M)(x-t)^{(M-I)}G(t) - C_m \frac{dG(t)}{dt} = 0 \qquad (3.33)$$

Eq. (3.32) has the solution on $G(t)$:

$$G(t) = K \exp\left(C_m^{-1} \int \Gamma^{-1}(M)(x-t)^{(M-I)} dt\right) \qquad (3.34)$$

Coming back to $F(x)$:

$$F(x) = \Gamma^{-1}(M) \int_0^x (x-t)^{(M-I)} G(t) dt \qquad (3.35)$$

Expanding eq. (3.32) in more terms:

$$C_m D^M F(x) + C_n D^N F(x) = F(x) \qquad (3.36)$$

Eq. (3.36) is similar to eq. (3.10), except that $\phi(x) = F(x)$. the same steps used to solve (3.10) apply to (3.36): consider $M = A + N$ and $F(x) = D^{-M}G(x)$:

$$C_m G(x) + C_n D^N D^{-A} G(x) = D^{-M} G(x) \qquad (3.37)$$

Opening (3.11) in terms of the integrals:

$$C_m G(x) + C_n \Gamma^{-1}(A) \int_0^x (x-t)^{(A-I)} G(t) dt = \Gamma^{-1}(M) \int_0^x (x-t)^{(M-I)} G(t) dt$$

$$(3.38)$$

Putting all together inside the integral, the equation obtained is:

$$\frac{dG(t)}{dt} + C_m^{-1}C_n\Gamma^{-1}(A)(x-t)^{(A-I)} - \Gamma^{-1}(M)(x-t)^{(M-I)} = 0 \qquad (3.39)$$

Eq. (3.39) is a variable separable equation whose solution is:

$$G(t) = K\exp\left(-\int C_m^{-1}C_n\Gamma^{-1}(A)(x-t)^{(A-I)} - \Gamma^{-1}(M)(x-t)^{(M-I)}dt\right)$$
$$(3.40)$$

The step 5 shows us the final solution:

$$F(x) = \Gamma^{-1}(M)\int_0^x (x-t)^{(M-I)}G(t)dt \qquad (3.41)$$

Generalizing eq. (3.36) to:

$$\sum_{n=1}^n C_{m_n}D^{M_n}F(x) = F(x) \qquad (3.42a)$$

$$(C_{m_1}D^{M_1} + C_{m_2}D^{M_2} + \cdots + C_{m_n}D^{M_n})F(x) = F(x) \qquad (3.42b)$$

The same procedures used to solve eq. (3.21), based on (3.22) also applies:

$$(C_{m_1} + C_{m_2}D^{-A_1} + \cdots + C_{m_n}D^{-A_{n-1}})G(x) = D^{-M_1}G(x) \qquad (3.43)$$

Following steps 3 – 5, we obtain:

$$\frac{dG}{dt} + C_{m_1}^{-1}\left[\Gamma^{-1}(M_1)(x-t)^{(M_1-I)} + \sum_{n=2}^n \Gamma^{-1}(A_{n-1})C_{m_n}(x-t)^{(A_{n-1}-I)}\right]G(t) = 0$$
$$(3.4)$$

Eq. (3.26) is a separable equation with solution:

$$G(t) = K\exp\left[-\int C_m^{-1}\left(\Gamma^{-1}(M_1)(x-t)^{(M_1-I)} + \sum_{n=2}^n \Gamma^{-1}(A_{n-1})C_{m_n}(x-t)^{(A_{n-1}-I)}\right)dt\right] \qquad (3.45)$$

Finally, we have $F(x)$:

$$F(x) = \Gamma^{-1}(M_1)\int_0^x (x-t)^{(M_1-I)}G(t) \qquad (3.46)$$

3.1.2. Homogeneous equations

Making $\phi(x) = 0$ in eq. (3.8):

$$C_m D^M F(x) + C_n D^N F(x) = 0 \qquad (3.47)$$

Following the same previous procedures, that is, considering $F(x) = D^{-M}G(x)$ and using properties (2.15) with $M = A + N$:

$$C_m G(x) + C_n D^N D^{-A-N} G(x) = 0 \qquad (3.48)$$

Eq. (3.48) is identical to eq. (3.11), except that $\phi(x) = 0$. Applying the same procedures used to solve eq. (3.11), the resulting differential equation is:

$$\frac{dG(t)}{dt} + C_m^{-1} C_n \Gamma^{-1}(A)(x-t)^{(A-I)} G(t) = 0 \qquad (3.49)$$

$$G(t) = K\exp\left(-\int C_m^{-1} C_n \Gamma^{-1}(A)(x-t)^{(A-I)} dt\right) \qquad (3.50)$$

Returning to $F(x)$:

$$F(x) = \Gamma^{-1}(M) \int_0^x (x-t)^{(M-I)} G(t)\, dt \qquad (3.51)$$

with $A = M - N$. Generalizing (3.48) for multiple matrix orders:

$$\sum_{n=1}^{n} C_{m_n} D^{M_n} F(x) = 0 \qquad (3.52a)$$

$$(C_{m_1} D^{M_1} + C_{m_2} D^{M_2} + \cdots + C_{m_n} D^{M_n}) F(x) = 0 \qquad (3.52b)$$

The process of solution of eqs. (3.52) is identical to eqs. (3.21), but setting up $F(x) = D^{-M_1} G(x)$ and explicating the order $M_1$ in terms of subsequent matrices $A_n$: $M_1 = M_n + A_{n-1}$, eqs. (3.22). The separable differential equation is:

$$\frac{dG(t)}{dt} + C_{m_1}^{-1} \sum_{n=2}^{n} \left( C_{m_n} \Gamma^{-1}(A_{n-1})(x-t)^{(A_{n-1}-I)} G(t) dt \right) = 0 \qquad (3.53)$$

The solution of eq. (3.53) is:

$$G(t) = K \exp\left(-\int C_m^{-1} \sum_{n=2}^{n} C_n (x-t)^{(A_{n-1}-I)} \Gamma^{-1}(A_{n-1}) dt\right) \qquad (3.54)$$

Finally, for $F(x)$:

$$F(x) = \Gamma^{-1}(M_1) \int_0^x (x-t)^{(M_1-I)} G(t) dt \qquad (3.55)$$

### 3.1.3 Equations with iterated orders

Consider the following equation:

$$D^M D^N F(x) = \phi(x) \tag{3.56}$$

With $\phi(x)$ a known continuous function, and $M, N$ commutable matrices. Eq. (3.56) is a type of matrix order differential equation with iterated matrix order, that is, there is a matrix order derivative of a matrix order derivative of a function, in a similar fashion to ordinary consecutive higher order derivatives, $\frac{d}{dx}\left(\frac{d}{dx}\right)$. Let us assume the solution $F(x)$ is the matrix integral of order $N$ of another function $G(x)$, by following the same procedures previous applied to solve the remaining equations:

$$F(x) = D^{-N} G(x) \tag{3.57}$$

Inserting (3.57) into (3.56) and using property (2.15):

$$D^M D^N D^{-N} G(x) = D^M G(x) = \phi(x) \tag{3.58}$$

Eq. (3.58) is identical to eq. (3.3). The solution is:

$$G(x) = \Gamma^{-1}(M) \int_0^x (x-t)^{(M-I)} \phi(t)\, dt \tag{3.59}$$

Now, retrieving $F(x)$ from (3.57):

$$F(x) = \Gamma^{-1}(N) \int_0^x (x-\tau)^{(N-I)} G(\tau)\, d\tau \tag{3.60}$$

Eq. (3.60) uses the variable $\tau$ only to not cause confusing with the previous integration variable $t$. In fact, if writing eq. (60) at full:

$$F(x) = \Gamma^{-1}(N) \int_0^x (x-\tau)^{(N-I)} \left[ \Gamma^{-1}(M) \int_0^\tau (\tau-t)^{(M-I)} \phi(t) \right) dt \bigg] d\tau \tag{3.61}$$

This same procedure of replacing the solutions for $G(t)$ inside the integral equations for $F(x)$ is possible in the preceding equations. However, for some cases the expression size may harm its reading.

Generalizing (3.56) to:

$$D^{M_n} D^{M_{n-1}} \dots D^{M_2} D^{M_1} F(x) = \phi(x)$$
$$\tag{3.62}$$

with all the $M_n$ as commutable matrices among them. Making:

$$F = D^{-M_1} G_1(x) \tag{3.63}$$

and inserting it in eq. (3.62):

$$D^{M_n} D^{M_{n-1}} \ldots D^{M_2} D^{M_1} D^{-M_1} G_1(x) = \phi(x) \tag{3.64}$$

$$D^{M_n} D^{M_{n-1}} \ldots D^{M_2} G_1(x) = \phi(x) \tag{3.65}$$

Eq. (3.65) is also written as:

$$\prod_{n=k}^{1} D^{M_n} F(x) = \phi(x) \tag{3.66}$$

Subsequently, assuming:

$$G_1(x) = D^{-M_2} G_2(x) \tag{3.67a}$$

$$G_{n-1}(x) = D^{-M_n} G_n(x) \tag{3.67b}$$

Rearranging up to the *n*th order of matrix derivation (3.65):

$$D^{M_n} G_{n-1}(x) = \phi(x) \tag{3.68}$$

A clearer example for $n = 3$ is given below to solve eq. (3.68):

$$D^{M_3} D^{M_2} D^{M_1} F(x) = \phi(x) \tag{3.69}$$

Making:

$$\begin{aligned} F(x) &= D^{-M_1} G_1(x) \\ G_1(x) &= D^{-M_2} G_2(x) \end{aligned} \tag{3.70}$$

$$\begin{aligned} D^{M_3} D^{M_2} D^{M_1} D^{-M_1} G_1(x) &= D^{M_3} D^{M_2} G_1(x) = \phi(x) \\ D^{M_3} D^{M_2} D^{-M_2} G_2(x) &= D^{M_3} G_2(x) = \phi(x) \end{aligned} \tag{3.71}$$

The consecutive solutions for eq. (3.71) are:

$$G_2(x) = \Gamma^{-1}(M_3) \int_0^x (x-t)^{(M_3-1)} \phi(t)\, dt$$
$$G_1(x) = \Gamma^{-1}(M_2) \int_0^x (x-t)^{(M_2-1)} G_2(t)\, dt \tag{3.72}$$

Finally, by applying (3.70):
$$F(x) = \Gamma^{-1}(M_1) \int_0^x (x-t)^{(M_1-1)} G_1(t)\, dt \tag{3.73}$$

If writing eq. (3.73) in a complete form, we have:

$$F(x) = \Gamma^{-1}(M_1) \int_0^x (x-t)^{(M_1-1)} \{\Gamma^{-1}(M_2)\Gamma^{-1}(M_3) \int_0^\tau (\tau - v)^{(M_3-1)} \phi(v)\, dv d\tau\}\, dt \tag{3.74}$$

Now, generalizing the solution (3.74) for the generalized eqs. (3.66) and (3.68), from $x_0$ to $x_n$ variables, we have:

$$F(x) = \prod_{n=0}^{k} \Gamma^{-1}(M_{n+1}) \int_0^{x_n} (x_n - x_{n+1})^{(M_{n+1}-1)} \phi(x_{n+1}) dx_{n+1} \tag{3.75}$$

Now, consider that $\phi(x)$ in eq. (3.62) – (3.66) is equal to the unknown function $F(x)$:
$$D^{M_n} D^{M_{n-1}} \ldots D^{M_2} D^{M_1} F(x) = F(x)$$
$$\prod_{n=k}^{1} D^{M_n} F(x) = F(x) \tag{3.66}$$

Solving an example similar to eq. (3.69):
$$D^{M_3} D^{M_2} D^{M_1} F(x) = F(x) \tag{3.69}$$

Making $F(x) = D^{-M_1} G_1(x)$, $G_1(x) = D^{-M_2} G_2(x)$ and using properties (2.12) and (2.15):
$$D^{M_3} G_2(x) = D^{-(M_1+M_2)} G_2(x) \tag{3.70}$$

Additionally, assuming there exists a third function $G_3$ such that $G_2(x) = D^{-M_3} G_3(x)$:
$$G_3(x) = D^{-(M_1+M_2+M_3)} G_3(x) = D^{-L} G_3(x) \tag{3.71}$$

Eq. (3.71) is similar to eq. (3.32) or eq.(3.6) with $\phi(x) = F(x)$. However, eq. (3.71) is a matrix integral equation, whereas (3.6) and (3.32) are matrix differential equations. Expanding (3.71) in terms of the definition (1.8):

$$G_3(x) = \Gamma^{-1}(L) \int_0^x (x-t)^{(L-I)} G_3(t)\, dt \tag{3.72}$$

Bringing $G_3(x)$ inside the integral as $\frac{dG(t)}{dt}$ provided that $G(0) = 0$ gives a separable differential equation:

$$\frac{dG_3(t)}{dt} - \Gamma^{-1}(L)(x-t)^{(L-I)} G_3(t) = 0 \tag{3.73}$$

The general solution for (3.73) is:

$$G_3(t) = K \exp\left(\int_0^x \Gamma^{-1}(L)(x-t)^{(L-I)} dt\right) \tag{3.74}$$

Back substituting $G_3, G_2, G_1$ and , we have:

$$G_2(x) = \Gamma^{-1}(M_3) \int_0^x (x-t)^{(M_3-I)} G_3(t)\, dt \tag{3.75}$$
$$G_1(x) = \Gamma^{-1}(M_2) \int_0^x (x-t)^{(M_2-I)} G_2(t)\, dt \tag{3.76}$$
$$F(x) = \Gamma^{-1}(M_1) \int_0^x (x-t)^{(M_1-I)} G_1(t)\, dt \tag{3.77}$$

$$F(x) = \Gamma^{-1}(M_1) \int_0^x (x-t)^{(M_1-I)} \left\{ \Gamma^{-1}(M_2) \int_0^t (t-\tau)^{(M_2-I)} \left[ \Gamma^{-1}(M_3) \int_0^\tau (\tau - v)^{(M_3-I)} \left( K \exp\left(\int_0^v \Gamma^{-1}(L)(v-u)^{(L-I)} du\right)\right) dv \right] d\tau \right\} dt \tag{3.77}$$

With $= M_1 + M_2 + M_3$. Generalizing solution (3.77) applied to eq.(3.66):

$$F(x) = \Gamma^{-1}\left(\sum_{j=1}^n M_j\right) \left[\prod_{n=0}^{k-1} \Gamma^{-1}(M_n) \int_0^{x_n} (x_n - x_{n+1})^{(M_{n+1}-I)} dx_n\right] \int_0^x (x_{n+1} - x_k)^{-\sum_{j=1}^n M_j} dx_k \tag{3.75}$$

with $L = \sum_{j=1}^{n} M_j$, provided that all the $M_n$ matrices be commutable.

3.2 Matrix order partial differential equations (MOPDE)

3.2.1 Equal order of the partial derivatives

Consider the following equation:

$$\frac{\partial^M}{\partial x^M} F(x, y) = \frac{\partial^M}{\partial y^M} F(x, y) \tag{3.76}$$

or, writing in a more concise notation:

$$\partial_x^M F(x, y) = \partial_y^M F(x, y) \tag{3.77}$$

This equation can be considered as a matrix order partial differential equation with equal derivative order. A possible way of solution is the use of the known method of variables separation, by assuming that:

$$F(x, y) = \Psi(x)\Psi(y) \tag{3.78}$$

Inserting this expression in the original equation, one has:

$$\Psi(y) \frac{\partial^M}{\partial x^M} \Psi(x) = \Psi(x) \frac{\partial^M}{\partial y^M} \Psi(y) = \kappa \tag{3.79}$$

Rearranging, we have:

$$\Psi^{-1}(x) \frac{\partial^M}{\partial x^M} \Psi(x) = \Psi^{-1}(y) \frac{\partial^M}{\partial y^M} \Psi(y) = \kappa \tag{3.80}$$

Now, if equating both sides of the equation with the constant matrix $\kappa$:

$$\frac{\partial^M}{\partial x^M} \Psi(x) = \kappa \Psi(x) \tag{3.81a}$$

$$\frac{\partial^M}{\partial y^M}\Psi(y) = \kappa\Psi(y) \qquad (3.81b)$$

Assuming there are functions $\eta(x)$ and $\zeta(y)$ in a similar fashion used to solve eq. (3.32):

$$\Psi(x) = D^{-M}\eta(x) \qquad (3.82a)$$
$$\Psi(y) = D^{-M}\zeta(y) \qquad (3.82b)$$

Replacing (3.82) into (3.81):

$$D^{-M}\eta(x) - \eta(x) = 0 \qquad (3.83a)$$
$$D^{-M}\zeta(y) - \zeta(y) = 0 \qquad (3.83b)$$

Opening eqs. (3.83) with the definition (1.8) and inserting $\eta(x)$ and $\zeta(y)$ inside the integral, since $\eta(t=0) = 0$ and $\zeta(\tau=0) = 0$ we have:

$$\frac{d\eta(t)}{dt} = \kappa(x-t)^{(M-I)}\Gamma^{-1}(M)\eta(t) \qquad (3.84a)$$
$$\frac{d\zeta(\tau)}{dt} = \kappa(y-\tau)^{(M-I)}\Gamma^{-1}(M)\zeta(\tau) \qquad (3.84b)$$

The general solution for (3.84) is:

$$\eta(t) = C_\eta \exp\left(\Gamma^{-1}(M)\int (x-t)^{(M-I)} dt\right) \qquad (3.85a)$$
$$\zeta(\tau) = C_\tau \exp\left(\Gamma^{-1}(M)\int (y-\tau)^{(M-I)} d\tau\right) \qquad (3.85b)$$

In eqs. (3.85) $C_\eta$ and $C_\eta$ are integration constant matrices. Applying (3.82) to recover $\Psi$, we find:

$$\Psi(x) = \kappa\Gamma^{-1}(M) \int_0^x (x-t)^{(M-I)}\eta(t)\, dt \qquad (3.86a)$$

$$\Psi(y) = \kappa\Gamma^{-1}(M) \int_0^y (y-\tau)^{(M-I)}\zeta(\tau)\, d\tau \qquad (3.86b)$$

Finally, using (3.78) the final solution for eq. (3.76) is:

$$F(x,y) = \kappa^2\Gamma^{-2}(M) \int_0^x \int_0^y (x-t)^{(M-I)}(y-\tau)^{(M-I)}\eta(t)\zeta(\tau)\, d\tau dt \qquad (3.87)$$

### 3.2.2 Different order of the partial derivatives

Consider the following partial equation with different matrix order:

$$\frac{\partial^M}{\partial x^M} F(x,y) = \frac{\partial^N}{\partial y^N} F(x,y) \qquad (3.88)$$

or, in a more concise notation:

$$\partial_x^M F(x,y) = \partial_y^N F(x,y) \qquad (3.89)$$

The method of variable separation also works to this problem as well as for eq. (3.76) Taking the same procedures in the aforementioned equal order case, we have:

$$\Psi^{-1}(x) \frac{\partial^M}{\partial x^M} \Psi(x) = \Psi^{-1}(y) \frac{\partial^N}{\partial y^N} \Psi(y) = \kappa \qquad (3.90)$$

$$\frac{\partial^M}{\partial x^M} \Psi(x) = \kappa \Psi(x) \qquad (3.91a)$$

$$\frac{\partial^N}{\partial y^N} \Psi(y) = \kappa \Psi(y) \qquad (3.91b)$$

The solution for $\Psi(x)$ is identical to eq. (3.86a). For $\Psi(y)$, assuming:

$$\Psi(y) = D^{-N} \zeta(y) \qquad (3.92)$$

Which takes to:

$$\frac{d\zeta(\tau)}{dt} = \kappa (y - \tau)^{(N-1)} \Gamma^{-1}(N) \zeta(\tau) \qquad (3.93)$$

The solution for $\zeta(\tau)$ is equal to (3.85b), except for replacing $M$ by $N$:

$$\zeta(\tau) = C_\tau \exp\left( \Gamma^{-1}(N) \int (y - \tau)^{(N-1)} d\tau \right) \qquad (3.94)$$

Regrouping solutions for $\eta(t)$ and $\zeta(\tau)$ to find $\Psi(x)$ and $\Psi(y)$, respectively, and knowing that $\Psi(x)$ is the same solution found in eq. (3.86b):

$$\Psi(y) = \kappa\Gamma^{-1}(N) \int_0^y (y-\tau)^{(N-1)} \zeta(\tau)\, d\tau \tag{3.95}$$

Lastly, applying (3.78), the final solution for eq. (3.88) is:

$$F(x,y) = \kappa^2 \Gamma^{-2}(M) \int_0^x \int_0^y (x-t)^{(M-1)}(y-\tau)^{(N-1)} \eta(t)\zeta(\tau)\, d\tau dt \tag{3.96}$$

4. System of linear matrix order differential equations

Let us suppose as an example, a system of a linear matrix order differential equations as follows:

$$\begin{cases} D^N F(x) + G(x) = 0 \\ D^M G(x) - F(x) = 0 \end{cases} \tag{3.97}$$

The solution for (3.97) consist of finding functions $(x)$, $G(x)$ that verify the equations in (3.97). We can begin the solution by different manners. A possible starting point is to enunciate that $F(x)$ is the *N*-integral of some function $R(x)$ with $R(x) \neq 0$:

$$F(x) = D^{-N} R(x) \tag{3.98}$$

Inserting (3.98) in (3.97):

$$\begin{cases} R(x) + G(x) = 0 \\ D^M G(x) - D^{-N} R(x) = 0 \end{cases} \tag{3.99}$$

Whence we find $G(x) = -R(x)$. Again, Now considering also $G(x)$ depends of some function $H(x)$:

$$G(x) = D^{-M} H(x) \tag{3.100}$$

$$\begin{cases} R(x) + G(x) = 0 \\ H(x) - D^{-N} R(x) = H(x) + (D^{-N} D^{-M} H(x)) = 0 \end{cases} \tag{3.101}$$

And using property (2.12):

$$\begin{cases} R(x) + G(x) = 0 \\ H(x) + (D^{-N-M}H(x)) = 0 \end{cases} \quad (3.101)$$

The equation on $H(x)$ in (3.101) is written as:

$$\int_0^x \left( \Gamma^{-1}(N+M)(x-t)^{(N+M-)}H(t) + \frac{dH(t)}{dt} \right) dt = 0 \quad (3.102)$$

As the integral is zero (remembering that $H(0) = 0$), the term inside the integral is also zero and it is a linear differential equation with separable variables. This solution is:

$$H(t) = K \exp\left(-\int (x-t)^{(N+M-I)} \Gamma^{-1}(N+M) dt\right) \quad (3.103)$$

From (3.100), we have for $G(x)$:

$$G(x) = \Gamma^{-1}(M) \int_0^x (x-t)^{(M-I)} H(t)\, dt = -R(x) \quad (3.104)$$

Consequently, for $F(x)$ as (3.98):

$$F(x) = -\Gamma^{-1}(N) \int_0^x (x-t)^{(N-I)} R(t) dt = \Gamma^{-1}(N) \int_0^x (x-t)^{(N-I)} G(t) dt \quad (3.105)$$

As any other system of equations, many other ways to work with this type of system enable to reach out faster and easier the final solution.

5. Conclusions

This work presented an extension of the fractional calculus by considering the concept of using the derivative order as a matrix. The principal properties derived from the fractional calculus were extended and proved for this matrix order calculus. This paper gives a special emphasis on the matrix order version of the Riemann-Liouville differential and integral operators. However, this approach is not restricted to this operator. It is possible to extend for other matrix order differintegrators based on their original fractional versions listed in Table 1. New investigations are necessary to check the validity of the theorems presented here for other operators. This work also presented exact results for some main types of matrix order differential equations. Future works on this topic will explore numerical solutions for those matrix order differential equations by extending classical numerical approximations to matrix order differintagration problems. Some possible applications in different areas are under investigation.

Acknowledgments

This work has not received any funding in their development and there is not any conflict of interest by the author.

Funding: This research did not receive any specific grant from funding agencies in the public, commercial, or not-for-profit sectors.